\documentclass[12pt]{article}
\usepackage{amsfonts}
\usepackage{mathrsfs}
\usepackage{amsmath,amssymb}

\baselineskip 5.5pt \lineskip 5.5pt \numberwithin{equation}{section}

\def\QED{\hfill$\Box$\par}
\def\cl{\centerline}
\def\vs{\vspace*}
\def\ni{\noindent}
\def\C{\mathbb{C}}
\def\Z{\mathbb{Z}}
\def\la{\lambda}

\def\bo{{\bf 1}}
\def\rar{\rightarrow}

\def\uu{\mathcal{U}}
\def\pp{\mathcal{P}}

\def\p-{\mathcal{P}_{\leq 0}}
\def\G{\mathfrak{g}}

\def\n{\mathfrak{n}}
\def\h{\mathfrak{h}}
\def\b{\mathfrak{b}}

\def\a{\alpha}

\begin{document}
\cl{{\bf\large WHITTAKER MODULES FOR GRADED LIE ALGEBRAS }}\vs{10pt}

\cl{Bin Wang}

\begin{abstract}{\footnotesize In this paper, we study Whittaker modules
for graded Lie algebras. We define Whittaker modules for a class of
graded Lie algebras and obtain a bijective correspondence between
the set of isomorphism classes of Whittaker modules and the set of
ideals of a polynomial ring, parallel to a result from the classical
setting and the case of the Virasoro algebra. As a consequence of
this, we obtain a classification of simple Whittaker modules for
such algebras. Also, we discuss some concrete algebras as examples.}
\end{abstract}


\noindent{\it{MR(2000) Subject Classification}: 17B10, 17B65,
17B68.}\vs{15pt}

\cl{\bf\S1.Introduction}
Let $\G$ be a Lie algebra that admits a decomposition $$\G = \b_-
\oplus \n $$where $\b_-, \n$ are two Lie subalgebras. Let $$\varphi
: \n \rar \C$$be a homomorphism of Lie algebras. For a $\G$ module
$V$ and $v\in V$, one says that $v$ is a Whittaker vector of type
$\varphi$ if $\n$ acts on $v$ through $\varphi$. A Whittaker module
is then defined to be a module generated by a Whittaker vector.

The category of Whittaker modules for a given algebra (say, $\G$)
admits an initial object and we call it a universal Whittaker module
(say, $M$). Then $M$ is isomorphic to $\uu (\b_-)$ as $\b_-$
modules. Let $Z$ stand for the center of $\G$, and $Z'' = Z \cap
\b_-$. Suppose $\G$ possesses the following properties:

1) for each ideal $I$ of $S(Z'')$, every whittaker vector of the
Whittaker module $M/IM$ is of form $p\bar{w}$, with $p \in S(Z'')$,
where $w$ is a Whittaker generator of M;

2) for each ideal $I$ of $S(Z'')$, every nontrivial submodule of
$M/IM$ admits a nonzero Whittaker vector.\vspace*{5pt}

\ni Then it is not hard to set up a correspondence between the set
of isomorphism classes of Whittaker modules and the one of ideals of
$S(Z'')$.

In this paper, we consider a class of infinite-dimensional Lie
algebras that has a grading over an ordered free abelian group $Q$.
Assume there is a homomorphism $\pi : Q \rar \Z$, preserving the
ordering. Then we define Whittaker modules for $\G$ through the
following decomposition
$$\mbox{$\G = (\bigoplus\limits_{\pi(\a) \leq 0}\G_\a) \bigoplus (\bigoplus\limits_{\pi(\a) >
0}\G_\a)$}$$ where the second part stands for $\n$. We assume
$dim(\G _\a ) \leq 1, \, \pi(\a)>0$ and $dim(\G/Z)_\a \leq 1,
\pi(\a) \leq 0$. Fix a nonzero element $x_\a$ of $\G_\a$, for each $
\a \in K\setminus R''$, such that $x_\a \notin Z''$.

Furthermore, assume that the set $\{ \a  \,|\, \pi(\a) > 0, \varphi
(x_\a) \neq 0 \}$ admits a unique maximal element, denoted by
$\a_{\varphi}$, and that for any $x_{\beta} \in \G_\beta , x_\beta
\notin Z, \pi(\beta)\leq 0 $, there exists a $y \in \n_{(\a_\varphi
-\beta)}$ such that $\varphi ([y, x_{\beta}])\neq 0$. Our main
result is to show that such a $\G$ satisfies the above properties 1)
and 2) (see \S3), and hence we can obtain a bijective correspondence
between Whittaker modules and ideals of $S(Z'')$ (see \S4). We also
study two examples, the Virasoro algebra and a type of $W$-algebra.

Whittaker modules were first discovered for $\mathfrak{sl}_2{(\C )}$
by Arnal and Pinzcon in \cite{AP}. Block showed, in \cite{B}, that
the simple modules for $\mathfrak{sl}_2(\C )$ consist of highest
(lowest) weight modules, Whittaker modules and a third family
obtained by localization. This illustrates the prominent role played
by Whittaker modules.

Kostant defined Whittaker modules for an arbitrary
finite-dimensional complex semisimple Lie algebra $\mathfrak{g}$ in
\cite{K}, and showed that these modules, up to isomorphism, are in
bijective correspondence with ideals of the center
$Z(\mathfrak{g})$. In particular, irreducible Whittaker modules
correspond to maximal ideals of $Z(\mathfrak{g})$. In the quantum
setting, Whittaker modules have been studied by Sevoystanov for
$\uu_h(\mathfrak{g})$ \cite{S} and by M. Ondrus for
$U_q(\mathfrak{sl}_2)$ in \cite{O}. Recently Whittaker modules have
also been studied by M. Ondrus and E. Wiesner for the Virasoro
algebra in \cite{OW}, X. Zhang and S. Tan for
Schr\"{o}dinger-Virasoro algebra in \cite{ZT}, K. Christodoulopoulou
for Heisenberg algebras in \cite{C}, and by G. Benkart and M. Ondrus
for generalized Weyl algebras in \cite{BO}.

 The paper is organized in the following way. In section 2,
we define Whittaker vectors and Whittaker modules for a class of Lie
algebras, and also construct a universal Whittaker module for them.
Then the whittaker vectors in a Whittaker module are examined in
section 3 and the irreducible Whittaker modules are classified in
section 4. In the last section we discuss some examples.

\vs{18pt}

\cl{\bf\S2. Preliminaries}

\ni 2.1. {\bf Q-graded Lie algebras} \vspace{6pt}

\ni 2.1.1. Let $V$ be a vector space over $\C$ and $Q$ a free
abelian group. By a $Q$-grading of $V$ we will understand a family
$\{V_\a | \a\in Q \}$ of subspaces of $V$ such that $V =
\oplus_{\a\in Q}V_{\a}$. For a nonzero vector $v \in V_\a$, we say
$v$ is a homogeneous vector of degree $\a$. Let $\G$ be a Lie
algebra over $\C$ and let $\{ \G_\a \,|\, \a\in Q \}$ be a grading
of $\G$ (as a vector space). Call $\G$ a $Q$-graded Lie algebra if
$[ \G_\a, \G_\beta ] \subset \G_{\a +\beta}$, for all $\a, \beta \in
Q$.

Now suppose $Q$ is totally ordered abelian group by the ordering,
$\leq$,  which is compatible with its group structure. Given a
$Q$-graded Lie algebra, $\G = \oplus_{\a\in Q}\G_{\a}$, and a
homomorphism of abelian groups $\pi : Q \rar Z$ that preserves the
ordering, write $\G_m = \sum\limits_{\pi (\a )=m} \G_{\a} $. Then
the Lie algebra $\G = \oplus_{i\in \Z}\G_i$ can be viewed as a
$\Z$-graded Lie algebra too. Set
\begin{eqnarray*}
&&\n = \mbox{$\bigoplus\limits_{m > 0}$}\G_m ,\\
&&\n_{-} = \mbox{$\bigoplus\limits_{m < 0}$}\G_m , \\
&&\b_- = \mbox{$\bigoplus\limits_{m\leq 0}$}\G_m ,\\
&&\mathfrak{h} = \G_0 .
\end{eqnarray*}
So $\G$ has a decomposition$$\G = \n_{-} \oplus \h \oplus \n = \b_-
\oplus \n .$$ Set $Q'=\{ \a \in Q \,|\, \pi (\a) >0 \}$ and $Q''= \{
\a \in Q \,|\, \pi (\a) \leq 0 \}$, then obviously $\n$ (resp.
$\b_-$) is $Q'$-graded (resp. $Q''$-graded ) Lie algebra.

\ni {\bf Definition}\ \  {\it  We call a $Q$-graded Lie algebra $\G$
a $Q$-good Lie algebra if
$$dim(\G _\a ) \leq 1, \, \forall \, \a \in Q' \,\, \mbox{and}\,\, \dim ((\G /Z)_\a ) \leq 1, \, \forall \, \a \in Q''  $$ where $Z$ is the center of $\G$.}

We assume, and will always assume for the rest of the paper, that
$\G$ is a $Q$-good Lie algebra. \vspace{6pt}

\ni 2.1.2.  We fix some notation here for the rest of the paper.
Write $\n_\a = \G_\a, \a \in Q'$. Let $K = \{\a \in Q \,|\, \G_\a
\neq 0 \}$. Set $K'=K\cap Q', K''=K\cap Q''$. The center of a graded
Lie algebra is always a graded ideal. So $Z$ is a graded ideal of
$\G$, that is, $Z = \bigoplus\limits_{\a \in Q}Z\cap \G_\a$. Let $R
= \{ \a \in K \,|\, Z\cap \G_\a = \G_\a \}$, and $R' = R\cap K', R''
= R\cap K''. $ Furthermore, define $Z''= Z\cap \b_- $. Fix a nonzero
element $x_\a$ of $\G_\a$, for each $ \a \in K\setminus R''$, such
that $x_\a \notin Z''$. \vspace{10pt}

\ni 2.2. {\bf{Partitions}}.\vspace{6pt}

\ni 2.2.1. Let $\Lambda$ be a totally ordered set. We define a
partition of $\Lambda$ to be a non-decreasing sequence of elements
of $\Lambda$,
$$\mu=(\mu_1,\mu_2,\cdots,\mu_r),\,\, \mu_1\leqslant\mu_2\leqslant\cdots\leqslant\mu_r.$$
Denote by $\pp (\Lambda )$ the set of all partitions. For $\lambda
=(\la_1,\cdots ,\la_r) \in \pp (\Lambda )$, we define the length of
$\la$ to be $r$, denoted by $\ell (\la)$, and  for $\a\in \Lambda$,
let $\la (\a)$ denote the number of times $\a$ appears in the
partition. Clearly any partition $\la$ is completely determined by
the values $\la (\a), \a\in \Lambda$. If all $\la (\a)=0$, call
$\la$ the null partition, denoted by $\bar{0}$. Note that $\bar{0}$
is the only partition of length $=0$. We consider $\bar{0}$ an
element of $\pp(\Lambda )$.

Back to the situation of $\G$. Define the symbols $x_{\la}$, for all
$\bar{0} \neq \la \in \pp(K\setminus R'')$, by
$$x_{\la}=x_{\la_1}x_{\la_2}\cdots x_{\la_r}=\prod\limits_{\a \in
K\setminus R''}x_\a^{\la(\a)} \in \uu(\G)$$ and set $x_{\bar{0}}=1
\in \uu{(\G )}$, where $\uu(\G)$ is the universal enveloping algebra
of $\G$. By PBW theorem, we know that $\{x_{\la} \,|\, \la \in
\pp(K\setminus R'') \}$
 form a basis of $\uu (\G )$ over $S(Z'')$.
\vspace{6pt}

\ni 2.2.2.\, $\uu{(\G)}$ (denoted by $\uu$) naturally inherits a
grading from the one of $\G$. Namely, for any $\a \in Q$, set
$\uu_\a =Span_{\C}\{x_1x_2\cdots x_k \,|\, x_i\in \G_{\a_i}, 1\leq
i\leq k, \sum\limits_{i=1}^{k}\a_i = \a \}$, and then $\uu =
\bigoplus\limits_{ \a \in Q} \uu_\a$ is a $Q$-graded algebra, i.e.
$\uu_\a \uu_\beta \subseteq \uu_{\a + \beta}$. Similarly, $\uu (\n
)$ (resp. $\uu ({\b_-})$) inherit a grading from $\n$, (resp. $\b_-$
). If $x \in \uu_\a$, then we say $x$ is a homogeneous element of
degree $\a$. Set $|\bar{0}|=0$ and $|\la |=\la_1 + \la_2 +\cdots +
\la_{{\ell}(\la)}, \, \forall \la\neq \bar{0}$. Then $x_\la$ is a
homogeneous element of degree $|\la |$. If $u (\neq 0 )$ is not
homogeneous but a sum of finitely many nonzero homogeneous elements,
then denote by $mindeg(u)$ the minimum degree of its nonzero
homogeneous components.

Now let us, for convenience, call any product of elements $x_\a^s$ (
$\a \in K\setminus R'', s\geq 0$) in $\uu$ and elements of $S(Z'')$
a monomial, of height equal to the sum of the various $s$'s
occurring. Then we have, by PBW theorem,\vspace{5pt}

\ni {\bf Lemma} \ \ {\it For $\a, \beta \in K\setminus R'',\,
t,k\geq 0$ , $x_{\beta}^t x_\a^k$ is a $S(Z'')$-linear combination
of $x_\a^k x_{\beta}^t$ along with other monomials of height $<
t+k$.} \QED

This allows us to make the following definition. If $ x \in \uu $ is
a sum of the monomials of height $\leq l$, we say $ht(x)\leq
l$.\vspace{6pt}

\ni 2.2.3. We need some more notation. For $\la =(\la_1,\la_2,\cdots
\la_r)\in \pp (K) , 0< i\leq r, 0\leq j < r, $write
\begin{eqnarray*}
&&\la \{i\}=(\la_1 ,\cdots , \la_i ), \, \la \{0\}=\bar{0} \\
&&\la [j]=(\la_{j+1},\cdots ,\la_r ), \, \la [r]=\bar{0}\\
&&\la <i>=(\la_1,\cdots ,\la_{i-1},\hat{\la_i},\la_{i+1},\cdots
\la_r )
\end{eqnarray*}\vspace{5pt}

\ni {\bf Lemma }\ \  {\it  Write $\uu''$ for $\uu (\b_-)$. Let
$0\neq x \in \n_\beta, 0\neq y \in \uu''_\gamma$ with $\pi (\beta)
> 0, \pi(\gamma) \leq 0$.

1) if $ s=\pi (\beta + \gamma ) > 0 $, then $[x, y] = \sum\limits_{
s\leq \pi(\a )\leq \pi (\beta) } u_\a $ with $u_\a =
\sum\limits_{i}v^{(\a,i)}w^{(\a,i)}$ where $w^{(\a,i)} \in \n_\a,
v^{(\a,i)} \in \uu'' _{\beta + \gamma - \a}$ and $ht(v^{(\a,i)}) <
ht(y)$ if $ v^{(\a,i)} \neq 0$;

2)if $ \pi (\beta + \gamma ) \leq 0 $, then $[x, y] = \sum\limits_{
0< \pi(\a )\leq \pi (\beta) } u_\a  + u$ with $u \in \uu''_{\beta +
\gamma}$ and $u_\a = \sum\limits_{i}v^{(\a,i)}w^{(\a,i)}$ where
$w^{(\a,i)} \in \n_\a, v^{(\a,i)} \in \uu'' _{\beta + \gamma - \a}$
and $ht(v^{(\a,i)}) < ht(y)$ if $ v^{(\a,i)} \neq 0$}.

\ni {\bf Proof}\ \  Write $y = \sum\limits_{\la}f_{\la}x_\la$ where
$f_{\la} \in S(Z'' ), \la \in \pp(K''\setminus R'')$.

But for any $x \in \n_\beta$, one has,
\begin{eqnarray*}
[x, x_\la ] &=& \sum_{i=1}^{\ell (\la )}x_{\la\{i-1\}}[x,
x_{\la_i}]x_{\la[i]}\\
&=& \sum_{i=1}^{\ell (\la )}x_{\la <i>}[x, x_{\la_i}] +
\sum_{i=1}^{\ell (\la )}x_{\la\{i-1\}} [[x, x_{\la_i}], x_{\la[i]}].
\end{eqnarray*}
Then one can easily deduce that the lemma follows. \QED
\vspace{10pt}

\ni 2.3 {\bf{Whittaker Module}} \ \  Let $\{x_{\a} \,|\, \a \in
K\setminus R'' \}$ be as in 2.1.2. \vspace{6pt}

\ni 2.3.1. {\bf{Definition}}. {\it  Given a homomorphism of Lie
algebras $\varphi: \n \rar \C$, for a $\G$-module $V$, a vector
$v\in V$ is said to be a Whittaker vector if $xv=\varphi(x)v$ for
all $x\in \n$. Furthermore, if $v$ generates $V$, we call $V$ a
Whittaker module of type $\varphi$ and $v$ a cyclic Whittaker vector
of $V$.}\vspace{6pt}

\ni 2.3.2. Given a homomorphism $\varphi:\n \rar\C$, define
$\C_\varphi$ to be the one-dimensional $\n$-module defined by the
action $xa= \varphi(x) a$ for all $x\in \n$ and $a\in\C$. Then the
induced $\G$-module
\begin{eqnarray*}
M_\varphi=\uu(\G)\otimes_{\uu(\n)}\C_\varphi,
\end{eqnarray*}
is a Whittaker module of type $\varphi$ with the cyclic Whittaker
vector $w=\bo\otimes 1$.  By PBW theorem, it's easy to see that
$\{x_{\la}\,|\,\la \in \pp (K''\setminus R'') \}$ is a basis of
$M_\varphi$ over $S(Z'' )$ .

Besides, for any  ideal $I$ of $S(Z'')$, define $L_{\varphi, I} =
M_\varphi/IM_\varphi$ and denote  by $p_I$ the canonical
homomorphism $M_{\varphi}\rar L_{\varphi, I}$. Then $L_{\varphi, I}$
is a Whittaker module for $\G$. The following lemma makes
$M_{\varphi}$ become a universal Whittaker module. \vspace{6pt}

\ni {\bf Lemma} {\it Fix $\varphi$ and $M_\varphi$ as above. Let $V$
be a Whittaker module of type $\varphi$ generated by a Whittaker
vector $w'$. Then there is a unique map $\phi: M_\varphi\rightarrow
V$ sending $w=1\otimes 1$ to $w'$.}

\ni{\bf Proof}.\ \ Uniqueness is obvious. Consider $u\in \uu (\G )$.
One can write, by PBW,
$$u=\sum\limits_\alpha b_\alpha n_\alpha,\,\, b_\alpha\in \uu (\b_- ), n_\alpha\in \uu(\n)$$
If $uw=0$, then $uw=\sum\limits_\alpha b_\alpha
\varphi(n_\alpha)w=0$, and therefore $\sum\limits_\alpha b_\alpha
\varphi(n_\alpha)=0$. Now it's easy to see that the map $\phi :
M_{\varphi} \rar V$, defined by $\phi(uw) = uw'$, is well
defined.\QED\vspace{6pt}

\ni 2.3.3. Let $A=S(Z'')$ and $I$ be an ideal of $A$. Write $M =
M_{\varphi}, \p- = \pp (K''\setminus R'')$ and $w' = p_I (w) \in
M/IM.$\vspace{6pt}

\ni {\bf Lemma} \ \ {\it $M/IM$ admits a basis, $\{ x_\la w'\,|\,
\la \in \p- \},$ over $A/I$.}

\ni{\bf Proof} \ \ Note that $M = \bigoplus\limits_{\la \in
\p-}Ax_\la w$. Hence,
$$M/IM = A/I \otimes_A M = A/I \otimes_A (\bigoplus_{\la \in
\p-}Ax_\la w)=\bigoplus_{\la \in \p-}(A/I) x_\la w' .$$ Then the
lemma follows immediately. \QED\vspace{6pt}

\ni 2.3.4. Assume now that $I$ is an ideal of $A = S(Z'')$, and $w'
= p_I(w)$. Write $\uu'$ for $\uu (\n)$, and $\uu''$ for $\uu
(\b_-)$. Then $V = M/IM$ has a natural grading as a vector space.
Namely, based on Lemma 2.3.2, let, for any $\a \in Q''$, $V_\a = \{
\sum\limits_{\la \in \p-}a_\la x_\la w' \,|\, a_\la \in (A/I), a_\la
= 0 \,\,\mbox{if}\,\,|\la | \neq \a \}$ and then clearly $V
=\bigoplus\limits_{\a\in Q''} V_{\a}$. We say that a nonzero
homogeneous vector $v$ in $M$ is of degree $\a$ if $v \in V_\a$. If
$v\, (\neq 0 )$ is not homogeneous but a sum of finitely many
nonzero homogeneous vectors, then define $mindeg(v)$ to be the
minimum degree of its nonzero homogeneous components. Meanwhile, for
any nonzero vector $v \in V$, let $d(v) = mindeg(v)$ and then there
uniquely exist $v_\a w' \in V_\a,  \a \geq d(v)$ such that $v =
\sum\limits_{d(v)\leq \a} v_{\a}w'$ where $v_\a = \sum\limits_{|\la
|=\a}a_{\la}x_\la$ for $\la \in \p-, a_\la \in A/I$ and $v_{d(v)}
\neq 0$. Then define $\ell (v) = ht( v_{d(v)})$. Note that $V$ can
also be equipped with a $\Z$-grading through $\pi : Q \rar \Z$. With
this grading, we can introduce notation $mindeg_1 (v)$ and $\ell_1
(v)$, for each $0\neq v \in V$, parallel to $mindeg(v)$ and $\ell
(v)$ respectively.

\vs{18pt}

\cl{\bf\S3. Whittaker Vectors}

Let $\G$ be a given $Q$-good Lie algebra. In this section, we
characterize the Whittaker vectors in a Whittaker module. Let $M =
M_{\varphi}$ where $\varphi$ is a fixed homomorphism of Lie algebras
from $\n \rar \C$, and let $w=\bo\otimes 1 \in M_{\varphi}$, and
$Z=\mathfrak{Z}(\G )$, the center of $\G$. Notation as in 2.1.2.
\vspace{10pt}

\ni 3.1.  {\bf Definition} \ \ {\it We say $(\G, \pi )$ is
nonsingular with respect to $\varphi$ if

1)\, the set $\{ \a \in K' \,|\, \varphi (x_\a) \neq 0 \}$ admits a
unique maximal element (denoted by $\a_{\varphi}$),

2)\, for every $x_{\beta}, \beta \in K''\setminus R''$, there exists
a $y \in \n_{(\a_\varphi -\beta)}$ such that $\varphi ([y,
x_{\beta}])\neq 0,$}

\ni where $\{ x_\a \,|\, \a \in K\setminus R'' \}$ is as in 2.1.2.
\vspace{10pt}

\ni 3.2. Assume that $I$ is an ideal of $A = S(Z'')$. Set $V = M/IM$
and $w' = p_I ( w )$. Write $\p- = \pp (K''\setminus R'')$. Then we
have the following lemma. \vspace{6pt}

\ni 3.2.1.  {\bf Lemma} \ \ {\it Suppose $(\G, \pi )$ is nonsingular
with respect to $\varphi$. Then every Whittaker vector of $V$ is of
form $pw'$ with $p \in A = S(Z'')$.}

\ni{\bf Proof} \ \ Suppose $w''$ is a Whittaker vector of $V$. We
can write, by Lemma 2.3.3,
$$w'' = \sum_{\la \in \p-}p_{\la}x_{\la}w'$$ where $ p_{\la} \in
A/I$. Suppose there is at least a $\bar{0}\neq \la \in \p-$ s.t.
$p_{\la}\neq 0$. Otherwise, nothing needs to show. Meanwhile,
without loss of generality, we may assume $p_{\bar{0}} = 0.$\,
Indeed, if $p_{\bar{0}}\neq 0$, we may replace $w''$ by
$w''-p_{\bar{0}}w'.$

let $T = min\{ |\la| \,|\, p_{\la}\neq 0 \}, A'=\{ \la \,|\, |\la| =
T, p_{\la}\neq 0 \}$. Set $l = max\{\ell(\la) \,|\, \la \in A' \}$
(clearly $l
> 0$), $B=\{ \la \,|\, \la \in A', \ell(\la)= l \}$, and
$\sigma_{0} = min\{\la_1 \,|\, \la \in B \}$. Furthermore, let
$B'=\{ \la \,|\, \la_1= \sigma_{0}, \la \in B \}$. Note that for any
$\la \in B$, we have $\la_i \geq \sigma_{0}$, and that if
$\la_i=\sigma_{0}$ for some $i, \la \in B$, then $\la \in B'$.

Let $\gamma = \a_{\varphi} - \sigma_{0}$ and let $y \in \n_\gamma$
be such that $\varphi ([y, x_{\sigma_{0}}])\neq 0$.  Note that
$\gamma + \la_i
> \a_\varphi$, for every $\la \in B$, unless $\la_i = \sigma_{0}$. Observe that
$$ yw''-\varphi(y)w'' = \sum_{\la\in \p-'}p_{\la}[y, x_{\la}]w'=D_1 + D_2$$
where $D_1=\sum\limits_{\la \in A'}p_{\la}[y,x_{\la}]w'$, and
$D_2=\sum\limits_{\la \notin A'}p_{\la}[y,x_{\la}]w'$. Since $uw'=0,
\forall u \in \n_\a$ with $\a > \a_{\varphi}$, we have $mindeg(D_2)
> T + \gamma - \a_{\varphi} =T-\sigma_{0}$, by Lemma 2.2.3, if $D_2 \neq 0$.

Case i), $l=1$.

Then $A'$ contains only one partition $\la$ with $\la = (\la_1 )$
and $\la_{1} = T $. So $A = B$ and $T=\sigma_0$. Now $D_1 = p_{\la}
\varphi([y, x_{\sigma_{0}}])w'\neq 0$ and it is of degree $0 (= T -
\sigma_0).$
 So, $yw''-\varphi(y)w'' = D_1 + D_2 \neq 0$ which contradicts with
the assumption that $w''$ is a Whittaker vector.

Case ii), $l > 1$.

Write $D_1 = D_1' + D_1''$ where
\begin{eqnarray*}
&&D_1' = \sum_{\la \in B}p_{\la}[y,
x_{\la}]w'\\
&&D_1'' = \sum_{\la \in A'\setminus B}p_{\la}[y, x_{\la}]w'
\end{eqnarray*}
Clearly either $mindeg(D_1'')>T-\sigma_{0}$ or $\ell(D_1'') < l-1$,
by Lemma 2.2.3, if $D_1''\neq 0$.

Meanwhile, we have
\begin{eqnarray*}
D_1'&=&\sum_{\la \in B}\sum_{i=1}^{\ell(\la
)}p_{\la}x_{\la\{i-1\}}[y,
x_{\la_i}]x_{\la[i]}w'\\
&=& E + F ,
\end{eqnarray*}
where \begin{eqnarray*} E&=&\sum_{\la \in B}\sum_{i=1}^{\ell(\la
)}p_{\la}x_{\la<i>}[y, x_{\la_i}]w'\\
F&=&\sum_{\la \in B}\sum_{i=1}^{\ell(\la )}p_{\la}x_{\la\{i-1\}}[[y,
x_{\la_i}], x_{\la[i]}]w'.
\end{eqnarray*}
It's easy to see that either $mindeg(F)>T-\sigma_{0}$ or
$\ell(D_1'') < l-1$, again by Lemma 2.2.3, if $ F\neq 0$.

However, on the other hand, since $[y, x_{\la_i}]w' = 0$ whenever
$\gamma + \la_i > \a_{\varphi}$, we have
\begin{eqnarray*}E &=& \sum_{\la \in
B'}\sum_{\la_i= \sigma_{0}} p_{\la}x_{\la<i>}[y, x_{\la_i}]w',\\
&=& \sum_{\la \in B'}\la ( \sigma_{0} )p_{\la}x_{\la<1>}\varphi([y,
x_{\sigma_{0}}])w'
\end{eqnarray*}
Clearly the terms occurring in the above sum are all nonzero
homogenous vectors of degree $= T - \sigma_{0}$ and of length
$=l-1$, and they are independent from each other. Hence $mindeg(E)=
T-\sigma_0, \ell(E) = l-1$. Therefore, $ yw''-\varphi(y)w'' = E + F
+ D_1'' + D_2 \neq 0$. This contradicts with the assumption that
$w''$ is a Whittaker vector. \QED \vspace{6pt}

\ni 3.2.2. {\bf Lemma} \ \ {\it Any nontrivial submodule of $V =
M/IM$ contains a nonzero Whittaker vector.}

\ni{\bf Proof} \ \  Let $V_1$ be a submodule of $V$. Suppose $V_1$
contains no nonzero Whittaker vector. Use the notation in 2.3.4. Let
$t = max\{ mindeg_1(v) \,|\, v\neq 0, v \in V_1 \}, l = min \{\ell_1
(v) \,|\, mindeg_1 (v) = t, v \in V_1, v\neq 0 \}$. Take a $u \in V$
such that $mindeg_1 (u)= t, \ell_1 (u) = l$ (clearly, $l > 0$ ).
Write $u = \sum\limits_{ 0\geq a \geq t}u_a $, where $ u_a =
\sum\limits_{\pi(|\la |) = a}p_{\la, a}x_\la w'$, with $p_{\la, a}
\in A/I$.  Since $u$ is not a Whittaker vector, there exists a $x
\in \n_\sigma$, for some $\pi(\sigma)
> 0$ such that $u' := xu-\varphi(x)u = \sum\limits_{0 \geq a \geq
t}[x, u_a]w' \neq 0$, where $[x, u_\a]$ stands for
$\sum\limits_{\pi(|\la |) = a}p_{\la, a}[x, x_\la]$. Note that $u'$
is contained in $V_1$. Then, it's easy to see $mindeg_1([x, u_\a]w')
\geq a \geq t$, by Lemma 2.2.3, if $[x, u_\a] \neq 0$. So we have
$mindeg_1(u') \geq t$ and hence $mindeg_1(u') = t $ for the
definition of $t$. In this case, $[x, u_{t}]w' \neq 0$ and
$mindeg_1([x, u_{t}]w') = t$. But this forces $ \ell_1 ([x, u_{t}]w'
) < ht(u_t) = \ell_1(u)$ (c.f. Lemma 2.2.3). Thus, $\ell_1 ( u' ) <
l $, which contradicts with the definition of $l$. \QED
\vspace{10pt}

\ni 3.3. {\bf Proposition} \ \ {\it Suppose $(\G, \pi )$ is
nonsingular with respect to $\varphi$. Let $N$ be a submodule of $M
= M_{\varphi}$. Then $N = IM$ for some ideal $I$ of $A = S(Z'')$.}

\ni{\bf Proof} \ \ Set $I = \{ x \in A \,|\, xw \in N \}$. One
immediately sees that $I$ is an ideal of $A$ and $IM \subseteq N$.
So we can view $N/IM$ as a submodule of $M/IM$. If $N\neq IM$, then
there exists $pw' \in N/IM$, with $pw' \neq 0, p \in A$, $(w' =
p_I(w))$ by Lemma 3.2.1 and 3.2.2. So $pw \in N$ and hence $p \in
I$. Therefore $pw \in IM$, which contradicts with the fact that $pw'
\neq 0$ in $N/IM$. Thus, $N = IM$.\QED \vs{18pt}

\cl{\bf\S4. \  Whittaker Modules}

Now we are ready to determine all the simple Whittaker modules and
set up a correspondence between the set of ideals of $S(Z'')$ and
the one of Whittaker modules. Notation as in 2.1.2.  Fix a
homomorphism $\varphi : \n \rar \C$, and let $M = M_{\varphi}, w =
\bo \otimes 1$. Let $A = S(Z'')$. \vspace{10pt}

\ni 4.1. {\bf Propostion} \ \ {\it Suppose $(\G, \pi )$ is
nonsingular with respect to $\varphi$. Then every nontrivial
submodule of a Whittaker module of type $\varphi$ contains a
nontrivial Whittaker submodule of type $\varphi$.}

\ni{\bf Proof} \ \ It follows immediately from Proposition 3.3 and
Lemma 3.2.2.\QED\vspace{10pt}

\ni4.2. The character $\varphi : \n \rar \C$ naturally extends to a
character of $\uu (\n )$. Let $\uu_{\varphi}(\n )$ be the kernel of
this extension so that $\uu (\n ) = \C \oplus \uu_{\varphi}(\n )$.
Hence,$$ \uu(\G) = \uu(\b_- )\otimes \uu(\n ) = \uu(\b_- ) \oplus
I_{\varphi}$$where $I_{\varphi}=\uu(\G )\uu_{\varphi}(\n )$. For any
$u \in \uu(\G )$, let $u^{\varphi} \in \uu(\b_- )$ be its component
in $\uu(\b_- )$ relative to the above decomposition of $\uu(\G )$.

If $V$ is a Wittaker module generated by a Whittaker vector $v$, let
$\uu_v(\G)$ (resp. $\uu_V(\G )$) be the annihilator of $v$ (resp.
$V$).Then we have, immediately, $V \simeq \uu(\G )/\uu_v(\G )$. Set
$A_V = A \cap \uu_V(\G )$.

 \ni 4.2.1. {\bf Proposition} \ \ {\it Suppose $(\G, \pi )$
is nonsingular with respect to $\varphi$. Let $V$ be any $\G$ module
that admits a cyclic Whittaker vector $v$. Then $$\uu_v(\G ) =
\uu(\G )A_V + I_\varphi .$$}

{\bf Proof}\ \ Obviously the right hand side of the equation is
contained in the left hand side. So it is enough to show the other
way around. Using the universal property of $M$, we can choose a
surjective homomorphism $\psi : M \rar V$ that sends $w=\bo \otimes
1$ to $v$. Let $Y = ker(\psi )$. Then $Y = IM$ for some ideal $I
\subseteq A$, by 3.3.

But for any $x \in A_V$, i.e. $xv=0$, one has $x^\varphi v=0$ and
hence $x^\varphi w \in Y$. Then $x^\varphi w = \sum_ip_ix_iw,\, x_i
\in \uu'', p_i \in I$. Thus, $x^\varphi = \sum\limits_ip_ix_i
\subset \uu I$. Clearly $I \subseteq A_V$, therefore $x^\varphi \in
\uu A_V$. So, $x \in \uu A_V + I_\varphi$. \QED \vspace{6pt}

\ni 4.2.2. {\bf Theorem} \ \ {\it Suppose $(\G, \pi )$ is
nonsingular with respect to $\varphi$. Then the correspondence $$V
\rar A_V$$sets up a bijection between the set of all the isomorphism
classes of Whittaker modules and the set of all the ideals of $A =
S(Z'')$.}

\ni{\bf Proof} \ \ Note that for any $I \subseteq A$, if let $V =
M/IM$, then $A_V = I$. Now the theorem follows
immediately.\vspace{6pt}

\ni 4.2.3. {\bf Corollary} \ \ {\it Suppose $(\G, \pi )$ is
nonsingular with respect to $\varphi$. For any maximal ideal
$\mathfrak{m} \in S(Z'')$, $L_{\varphi, \mathfrak{m}} =
M/\mathfrak{m}M$ is simple and every simple Whittaker module of type
$\varphi$ is of this form.}

\ni{\bf Proof} \ \ Observe that if $I \subseteq J \subseteq S(Z'')$,
then $M/JM$ is a quotient of $M/IM$. The corollary now follows from
Theorem 4.2.2 immediately.\QED  \vs{18pt}

\cl{\bf\S5. Applications}

We now discuss some examples.

\ni 5.1. {\bf \ Virasoro algebra}. The Virasoro algebra
$\mathcal{V}$ is an infinite-dimensional Lie algebra with a
$\C$-basis $\{\,L_n,\,c\,|\,n\in \Z\,\}$ and the following Lie
brackets:
\begin{eqnarray*}
&&[L_n,L_m\,]=(m-n)L_{m+n}+\frac{n^3-n}{12}\delta_{m+n,0}c,\\
&&[L_n, c]=0.
\end{eqnarray*}
and it has the following decomposition
\begin{eqnarray*}
\mathcal{V}=\n_+\oplus \h \oplus \n,
\end{eqnarray*}
where
\begin{eqnarray*}
&&\h=\mbox{Span}_\C\{L_0,c\},\\
&&\n=\mbox{Span}_\C\{L_n\,|\,n\in\mathbb{Z},n>0\},\\
&&\n_-=\mbox{Span}_\C\{L_n \,|\,n\in\mathbb{Z},n<0\}.
\end{eqnarray*}
A character $\varphi : \n \rar \C$ is said to be nonsingular if
$\varphi (L_1 )\neq 0, \varphi (L_2 )\neq 0.$

Let $Q = \Z$ and $\pi : Q \rar \Z$ be the identity map. Set
$\mathcal{V}_a = \C L_a \oplus \C\delta_{a,0}c,\, a\in \Z$. Clearly
$\mathcal{V}$ is a $Q$-good algebra and for every nonsingular
character $\varphi$ of $\n$, $(\mathcal{V}, \pi )$ is nonsingular
with respect to $\varphi$. Hence analogues to the results in \S 4
hold for the Virasoro algebra as well.\vspace{10pt}

\ni 5.2. Let $W' = W'(2,2)$ be a $W$-algebra that is an
infinite-dimensional Lie algebra with a $\C$-basis
$\{\,L_n,\,W_n,\,z, c \,|\,n\in \Z\,\}$ with $c$ and $z$ contained
in the center and the following Lie brackets:
\begin{eqnarray*}
&&[L_n,L_m\,]=(m-n)L_{m+n}+\frac{n^3-n}{12}\delta_{m+n,0}c,\\
&&[L_n,W_m]=(m-n)W_{m+n}+\frac{n^3-n}{12}\delta_{m+n,0}z,\\
&&[W_n, W_m ]=0.
\end{eqnarray*}
and it has the following decomposition
\begin{eqnarray*}
\mathcal{V}=\n_+\oplus \h \oplus \n,
\end{eqnarray*}
where
\begin{eqnarray*}
&&\h=\mbox{Span}_\C\{L_0, W_0, c, z \},\\
&&\n=\mbox{Span}_\C\{L_n, W_n, \,|\,n\in\mathbb{Z},n>0\},\\
&&\n_-=\mbox{Span}_\C\{L_n, W_n \,|\,n\in\mathbb{Z},n<0\}.
\end{eqnarray*}
A character $\varphi : \n \rar \C$ is said to be nonsingular if
$\varphi (L_i )\neq 0, \varphi (W_i )\neq 0, i=1,2.$

Let $Q = \Z \times \Z$ be equipped with the usual group structure
and an ordering given by
$$ (a,i) < (b,j) \,\, \mbox{if either}\,\, a < b
\,\,\mbox{or} \,\, a=b, i<j . $$Hence, $Q$ becomes an ordered group.

Set $W'_{(a,0)} = \C L_a \oplus \C \delta_{a,0} c, W'_{(a,1)} = \C
W_a \oplus \C \delta_{a,0} z$ and $W'_{(a,i)} = 0, i\neq 0,1$. Then
$W'$ becomes a $Q$-graded Lie algebra. Define $\pi : Q \rar \Z$ by
$\pi((a,i))=a$. Clearly $W'$ is a $Q$-good algebra and for every
nonsingular character $\varphi$ of $\n$, $(W', \pi )$ is nonsingular
with respect to $\varphi$. Hence analogues to the results in \S 4
hold for $W'(2,2)$. Pass to $\frac{W'}{(c-z)W'} = W(2.2)$ and one
gets the similar results for $W(2,2)$ that was treated in
\cite{WL}.

\vspace{6pt}


\begin{thebibliography}{9999}\vskip0pt\small
\parindent=2ex\parskip=-1pt\baselineskip=-1pt

\bibitem {AP} D. Arnal and G. Pinczon, On algebraically irreducible representations of he Lie algebra $\mathfrak{sl}_2$
{\it J. Math. Phys.} {\bf 15} (1974), 350--359.

\bibitem {BO} G. Benkart and M. Ondrus, Whittaker modules for
Generalized Weyl Algebras, arXiv:0803. 3570.

\bibitem {B} R. Block, The irreducible representations of the Lie
algebra $\mathfrak{sl}_2$ and of the weyl algebra, {\it Adv. Math.}
{\bf 39} (1981), 69--110.

\bibitem {C} K. Christodoulopoulou, Whittaker modules for Heisenberg algebras and imaginary
Whittaker modules for affine Lie algebras, {\it J. Alg.} {\bf 320}
(2008) 2871--2890.

\bibitem {K} B. Kostant, On Whittaker vectors and representation theory,
{\it Invent. Math.}, {\bf 48} (1978), 101--184.

\bibitem {O} M. Ondrus, Whittaker modules for $U_q({\mbox{sl}}_2)$, {\it J.
Alg.} {\bf 289} (2005), 192--213.

\bibitem {OW} M. Ondrus and E. Wiesenr, Whittaker modules for the Virasoro
algebra, arXiv:0805.2686.

\bibitem {S} A. Sevostyanov, Quantum deformation of Whittaker modules
and Toda lattice, {\it Duke Math. J.}, (2000), {\bf 204} 211--238.

\bibitem {WL} B. Wang, J. Li. Whittaker Modules for $W$-algebra $W(2,2)$ preprint,
arXiv:0902.1592.

\bibitem {ZT} X. Zhang and S. Tan, Whittaker modules and a class of new
modules similar as Whittaker modules for the
Schr\"{o}dinger-Virasoro algebra, arXiv:0812.3245v1.

\vs{10pt}

\ni {\footnotesize Department of Mathematics, Changshu Institute of
Technology, Changshu 215500, China, \, Email: binwang72@hotmail.com}




\end{thebibliography}
\end{document}